\documentclass[a4paper, 12pt]{article}
\usepackage{amsfonts}
\usepackage{amssymb}
\usepackage{amsmath}
\usepackage{amsthm}

\newtheorem{defi}{Definition}[section]
\newtheorem{proposition}[defi]{Proposition}
\newtheorem{theorem}[defi]{Theorem}
\newtheorem{lemma}[defi]{Lemma}

\newtheorem{corollary}[defi]{Corollary}
\newtheorem{remark}[defi]{Remark}

\newcommand{\la}{\longrightarrow}

\newcommand{\car}{\mbox{\rm char }}
\newcommand{\Aut}{\mbox{\rm Aut}}
\newcommand{\el}{\ell}
\newcommand{\Ad}{\mbox{\rm Ad}}
\newcommand{\Fr}{\mbox{\rm Fr}}
\newcommand{\Tr}{\mbox{\rm Tr}}
\newcommand{\Spec}{\mbox{\rm Spec}}
\newcommand{\Gl}{\mbox{\rm Gl}}
\newcommand{\Id}{\mbox{\rm Id}}

\newcommand{\F}{{\mathbb F}}
\newcommand{\N}{{\mathbb N}}
\newcommand{\Z}{{\mathbb Z}}
\newcommand{\PP}{{\mathbb P}}
\newcommand{\red}{\mbox{\rm red}}

\title{Automorphisms groups for $p$-cyclic covers of the affine line}
\author{Claus Lehr and Michel Matignon}
\begin{document}
\maketitle
\begin{abstract}
Let $k$ be an algebraically closed field of positive characteristic
$p>0$ and $C \to {\mathbb P}^1_k$ a $p$-cyclic cover of the projective line  ramified in
exactly one point. We are interested in the $p$-part of the full
automorphism group $\Aut_k C$. First we prove that these groups are
exactly the extraspecial $p$-groups and groups G which are subgroups
of an extraspecial group E such that $Z(E) \subseteq G$.
The paper also describes an efficient algorithm to compute the 
$p$-part of $\Aut_k C$ starting from an Artin-Schreier equation for
the cover $C \to {\mathbb P}^1_k$. 

The interest for these objects initially came from the study of the
stable reduction of $p$-cyclic covers over the $p$-adics.
There the covers $C \to {\mathbb P}^1_k$ naturally arise and their automorphism
groups play a major role in understanding the arithmetic monodromy.
Our methods rely on previous work by Stichtenoth \cite{St1}, \cite{St2}
whose approach we have adopted.

\end{abstract}

\section{Introduction}

When considering semi-stable models for $p$-cyclic covers of the projective line over a $p$-adic field $K$, we get as irreducible components of the special fiber
$p$-cyclic covers (we mean \'etale) of the affine line (see
\cite{Le1}, \cite{Ma}). An interesting arithmetico geometric object is
the monodromy, i.e.\ the minimal Galois extension $K'/K$ (say with
group $G$) such that a semi-stable model is defined over $K'$. A
classical result (see \cite{Deschamps}, \cite[theorem 4.44, p.551]{Liu}) asserts that this group acts faithfully on the special fiber of the semi-stable model as an automorphism group; so the complexity of the monodromy group is intimately related with the automorphism group of $p$-cyclic covers of the affine line over the residue field $k$.
The aim of this paper is first to report and then to complete the
literature on the subject 
(\cite{St1}, \cite{St2}, \cite{vdG-vdV}) with the objective to use the
results in order to study the monodromy group which occurs when
considering $p$-cyclic covers of the projective line over a $p$-adic
field (see \cite{Le-Ma1}, \cite{Le-Ma2}).

Our choice to disjoin this paper from the one on monodromy came from the new
interest in automorphism groups of curves 
(see \cite{Guralnick}, \cite{Leo}, \cite{Poonen1}, \cite{Poonen2}).

Let us discuss the content of the paper. In the following theorem we have 
gathered the bounds  proved in the paper.

\begin{theorem} \label{t1}
Let $k$ be an algebraically closed field of $\car p>0$
and $f(X)\in k[X]$ a polynomial of  
degree $m:=\deg f$ prime to $p$. Let $k(X,W)/k(X)$ be an extension
of degree $p$ defined by 
$W^p-W=f(X)$ and denote by $\infty$ the place in $k(X,W)$ above
$X=\infty$. We write
$\Aut k(X,W)$ for the full automorphism group, $\supset G_{\infty}(f)$ (the inertia
group at $\infty$) $\supset G_{\infty,1}(f)$ (the wild inertia group
at $\infty$). Let $g(f)=\frac{(m-1)(p-1)}{2}$ be the genus of $k(X,W)/k$ and
assume $g(f)>1$.
Then

a) {\rm \cite[Satz 4]{St2}, [Lemma \ref{44} i)]}  $|G_{\infty,1}(f)|\leq p(m-1)^2$.

b) {\rm \cite[Satz 4]{St2}, [Lemma \ref{44} ii)]}  If $(m-1,p)=1$, then $|G_{\infty,1}(f)|=p$.

c) {\rm \cite[Satz 5]{St2}, [Proposition \ref{415}]} If $f(X)=X^m$, then $|G_{\infty,1}(f)|=p$ for $m-1\neq p^s$
and $|G_{\infty,1}(f)|=p^{2s+1}$ for $m-1=p^s$.

d) {\rm [Lemma \ref{413}]} If $m-1=\el p^s,\ s>0,\ \el >1,\ (\el,p)=1 $, then $|G_{\infty,1}(f)|\leq p^{s+1}$ for $p>2$ and $|G_{\infty,1}(f)|\leq 2^s$ for $p=2$; 
moreover these bounds are optimal.

e) {\rm [Proposition \ref{417}]} If $\frac{|G_{\infty,1}(f)|}{g(f)}>\frac{p}{p-1}$ then $f(X)=XR(X)$, where $R(X)$ is an additive polynomial. 
\end{theorem}

In order to answer the question ``What is the structure of $G_{\infty,1}(f)$'',
we face the following: For a given $f\in Xk[X]$ as above determine the
set:

$$
S(f):=\{y\in k\ |\ \exists P(X,y)\in  Xk[X], \quad 
f(X+y)-f(X)-f(y)=(\Id -F)P(X,y)\}
$$

It is easy to see that giving $S(f)$ is equivalent to finding the additive polynomial 
$\Ad_f(Y)=\prod_{y\in S(f)}(Y-y)$. Here we give an algorithm
(cf. section \ref{algor}) which for a given
$f\in Xk[X]$ produces   the polynomial $\Ad_f(Y)$ and for $y\in S(f)$
the polynomial $P(X,y)\in Xk[X]$.
From this we deduce the structure of the group 
 $G_{\infty,1}(f)$. Moreover we use this algorithm in order to produce polynomials $f(X)\in k[X]$ or families $f_{\underline t}\in k[X,\underline t]$, with
prescribed group $G_{\infty,1}(f)$.
Next we describe the groups $G_{\infty,1}(f)$. 

Let us first define 3 sets of isomorphism classes of $p$-groups.

$$C_1:=\{ G\ |\ \exists N\  p{\mbox{\rm -cyclic and normal subgroup}}\subset G\ |\ 
G/N\ \mbox{\rm is } p \mbox{\rm -elementary abelian}\}$$

$$C_2:=\{ G\ |\ \exists E\ \mbox{\rm extraspecial group, } V\subset
\frac {E}{Z(E)}  \  \F_p \mbox{\rm -subspace}\ |\  G\simeq  \pi^{-1}(V),$$
$${\mbox{\rm where } }\ \pi:E\to \frac {E}{Z(E)}
\mbox{\rm is the canonical map}\}$$

$$C_3:=\{ G\ |\ \exists f\in Xk[X], \ (\deg f,p)=1 \mbox{ \rm such
that } G\simeq G_{\infty,1}(f)\}$$

\begin{theorem} The 3 classes $C_i,i=1,2,3$ are equal.
\end{theorem}

\section{Notations}

 Throughout this paper we use the following notation:

\begin{itemize} 

\item $k$ is an algebraically closed field of characteristic $p>0$.

\item $F$ design the Frobenius endomorphism for a $k$-algebra.
 
\item  $f(X)\in k[X]$ a polynomial, the equation $ W^p-W=f(X)$
defines an \'etale cover of the affine line that we denote $C_f$; moreover each \'etale cover of the affine line can be presented like this. Let $f(X)\in k[X]$; there is a unique polynomial 
$\red (f)(X)\in k[X]$ called the {\it reduced representative} of $f$ which is $p$-powers free, i.e.\
$\red (f)(X)\in k\bigoplus_{(i,p)=1}kX^i$ and such that $\red (f)(X)=f(X) \mod (F- \Id)
k[X]$.  Clearly the covers $C_f$ and $C_{\red (f)}$ are the same $p$-cyclic 
cover of the affine line. The curve $C_f$ is irreducible iff $\red (f)\neq 0$.  In the sequel we assume that $\red (f)\neq 0$, the degree of the reduced polynomial
$\red (f)$ is called the {\it conductor} of the cover it is prime to $p$ and equal to the degree of $f$ if $(\deg f,p)=1$.

\item  By $C_f$ we also denote the non singular projective curve with function field $k(X,W)$. If $m=\deg f$ is prime to $p$ then the genus of 
$C_f$ is $g(C_f)=\frac{(m-1)(p-1)}{2}$.

\item For $M\in C_f$ and $n\in \N$ let $L(nM)=\{\varphi \in k(C_f)\ |\ (\varphi)+nM\geq 0\}$, $\mathcal P(M):=\{ -v_M(\varphi)\ | \ \varphi\in \cup_{n\in \N}L(nM)-\{0\}\}$ is the {\it polar semi-group} at $M$.

\item We denote by $\infty \in C_f$ the point $X=\infty,\ W=\infty$ and by $M_{x,w}$ the point $X=x,\ W=w$ with $w^p-w=f(x)$.

\item We denote $G(f)=\Aut_k C_f$, $G_{\infty}(f)$ the inertia group at $\infty$ and
$G_{\infty,1}(f)$ the wild inertia group at $\infty$.
Let $\rho\in G(f)$ be such that
$\rho(X)=X$ and $\rho(W)=W+1$. Then $\rho$ generates a $p$-cyclic
subgroup in 
$G_{\infty}(f)$.
\end{itemize}

\section {Review and improvements of Stichtenoth's \\results.}

Stichtenoth proves the following:

\begin{theorem}[\cite{St2}, Satz 6]  \label{stich}
Assume the genus $g(C_f)\geq 2$, i.e.\ $\{m,p\}\neq \{2,3\}$. Then $\infty$ is the only point $M\in C_f$ such that $\mathcal P(M)=\deg f 
\N+p\N$, with the following two exceptions:

a) $m<p,\ m|1+p,\ f(X)=X^m$. In this case in addition to 
$\infty $, exactly the $p$ zeroes of $X$ have the same semi-group $\mathcal 
P(M)$. More precisely for $i\in \F_p$, let $\sigma_i$ be given by
$\sigma_i(X)=\frac{X}{(W-i)^{\delta}}$ and
$\sigma_i(W)=-\frac{1}{W-i}$, where $\delta=\frac{1+p}{m}$.
Then $\sigma_i\in \Aut_k C_f$
and $\sigma_i(M_{0,i})=\infty$, in particular $\Aut_kC_f$ acts transitively on the $p+1$ points $M$ such that $\mathcal P(M)=\mathcal P(\infty)$.

b) $f(X)=X^{1+p}$. Now exactly the $p^3$ points $M_{\alpha,\beta}$ where 
$\alpha^{p^2}=\alpha$ and $\beta^p-\beta=\alpha$ have the same semi-group at
$\infty$.
Let $\zeta\in k$ with $\zeta^{p^2-1}=-1$,

$$\sigma_{\alpha,\beta}(X)=\alpha +\frac{X}{\zeta^{1+p}W}, \quad
\mbox{\rm and} \quad 
\sigma_{\alpha,\beta}(W)=\beta +\frac {1}{\zeta^{2(1+p)}W}-\frac{\alpha^p}{\zeta^2}\frac{X}{W}.$$

Then $\sigma_{\alpha,\beta}\in \Aut C_f$ and $\sigma_{\alpha,\beta}(\infty)=
M_{\alpha,\beta}$. In particular $\Aut C_f$ acts transitively on the $p^2+1$
points $M$, such that $\mathcal P(M)=\mathcal P(\infty)$.

\end{theorem}

\begin{remark} \rm Stichtenoth considers equations of type $W^p+W=f(X)$,
the expressions of the automorphisms 
$\sigma_i$ and $\sigma_{\alpha,\beta}$ in case a) and  b) are then simpler. 
For a description of the full automorphism group in case b) we refer
to \cite{Leo}.
\end{remark}

We deduce the following 

\begin{proposition} \label{33} Let $f\in Xk[X]$, $g\in Tk[T]$ and $C_g: V^p-V=g(T)$ 
such that $\deg f,\ \deg g$ are prime to $p$ 
and $g(C_f)\geq 2$. We assume $\exists \sigma: C_f\to C_g$, a $k$-isomorphism.
Then $\exists \sigma \in \Aut C_f$ such that $\varphi\circ\sigma(\infty)=
\infty $ and $\varphi\circ\sigma$ descends to $\PP^1$, i.e.\ $\exists \tilde\varphi: 
\PP^1\to \PP^1$ and a commutative diagram 

$$
\begin{array}{ccc}
C_f & \stackrel{\varphi \circ \sigma}{\la} & C_g \\
\downarrow &  & \downarrow \\
\PP^1 & \stackrel{\tilde{\varphi}}{\la} & \PP^1.
\end{array}
$$

Moreover $\exists (a,b)\in (k^{\times},k)$ and $c\in \F_p^{\times}$ such
that $\tilde\varphi^\sharp (X)=aX+b$ and $cf(X)-g(aX+b)\in (F-\Id)k[X]$.
If $f$ is not in case a) or b) of Theorem \ref{stich} we can take
$\sigma =\Id $.
Reciprocally any triple $(a,b,c)\in k^{\times}\times k\times 
\F_p^{\times}$ such that $cf(X)-g(aX+b)\in (F-\Id)(k[X])$ induces such a 
commutative diagram.

\end{proposition}

\begin{proof} $\varphi^\sharp :k(C_g)\to k(C_f)$ is a $k$-isomorphism which induces a graded
isomorphism between the linear spaces $L(\infty \varphi(M))$ and $L(\infty M)$ and so the polar semi-groups are preserved. It follows that $\varphi^{-1}(\infty)
\in C_f$ has a polar semi-group equal to that of $\infty \in C_f$. We deduce 
the existence of $\sigma $ and can assume from now on that $f(\infty)=\infty$. 

Let $m=\deg f$, it is classical that $g(C_f)=\frac{(m-1)(p-1)}{2}$, so as 
$g(C_f)=g(C_g)$ it follows that $m=\deg f=\deg g$. Now we follow the
discussion in Stichtenoth's paper \cite[Satz 4]{St2}.

Case 1. $m>p$. Then $1,T$, resp. $1,T,...,T^{[\frac {m}{p}]},V$
is a $k$-basis for $L(p\infty)\subset k(C_g)$ resp. $L(m\infty)\subset
k(C_g)$. Idem  $1,X$ and $1,X,...,X^{[\frac {m}{p}]},W$ are a basis for
$L(p\infty)\subset k(C_f)$ resp. $L(m\infty)\subset k(C_f)$. As
$\varphi^\sharp (L(p\infty))=L(p\infty)$ resp. $\varphi^\sharp
(L(m\infty))=L(m\infty)$, we get $\exists a\in k^{\times},b\in k$,
such that $\varphi^\sharp (T)=aX+b$ and $\varphi^\sharp (V)=cW+Q(X)$ where $c\in k^{\times}$ and $Q(X)\in k[X]$ with $\deg Q(X)\leq  [\frac {m}{p}]$. 
Then $(c^p-c)W+c^pf(X)+(Q(X)^p-Q(X))=g(aX+b)$, so
$c^p=c$ and $cf(X)-g(aX+b)\in (F-\Id)k[X]$.
Reciprocally such relations define an isomorphism.

Case 2. $m<p$. It works similarly, namely $1,V$, resp. $1,V,...,V^r,T$ with 
$rm<p<(r+1)m$ is a basis for $L(m\infty)$ resp. $L(p\infty)$ in
$k(C_g)$. Then
$\varphi^\sharp (T)=aX+P(W)$, where $P(W)\in k[W]$ with $m\deg P<p$
and $\varphi^\sharp (V)=cW+Q$ where $c\in k^{\times}$ and $Q\in k$. We get the equation $(c^p-c)W+c^pf(X)+(Q(X)^p-Q(X))=g(aX+P(W))=a^mX^m+(ma^{m-1}P(W)+b_{m-1})X^{m-1}+...+P(W)^m$,
where $g(T)=T^m+b_{m-1}T^{m-1}+...$. Comparing the $W$ degrees it follows that
$P(W)=b\in k$ and $c^p=c$. Now we get again the condition $cf(X)-g(aX+b)\in (F-\Id)k[X]$.
\end{proof}

\begin{corollary} \label{34}  Let $f\in Xk[X]$ such that $g(C_f)\geq 1$ and $S(f):=\{ (a,b,c)\in k^{\times}\times k\times \F_p^{\times}\ | \ \exists P_{a,b,c}(X)\in Xk[X], \ 
cf(X)-f(aX+b)=P_{a,b,c}(X)^p-P_{a,b,c}(X)\}$ moreover such a polynomial $P_{a,b,c}$ in
$Xk[X]$ is uniquely determined by the relation. Then for $(a,b,c)\in S(f)$, the 
formulas $\sigma_{a,b,c}(X)=aX+b$, $\sigma_{a,b,c}(W)=cW+P_{a,b,c}(X)$ define
 an automorphism of $C_f$ which lies in $G_{\infty}(f)$. 
Further $G_{\infty}(f)=<\rho,\sigma_{a,b,c}>$ for $(a,b,c)\in S(f)$. 

\end{corollary}

\begin{proof}  If $g(C_f)>1$, this follows from the proposition.
If  $g(C_f)=1$ then $\{p,m\}=\{2,3\}$; if $p=2$ then $f(X)= aX+X^3$ and if $p=3$ then $f(X)=aX+X^2$; in each case we can give a basis for $L(m\infty)$ and 
$L(p\infty)$ of the same kind of those given in the two cases in the proof of 
the proposition; we conclude in the same way. Note that if $g(C_f)=0$ then 
$G_{\infty}(f)=k^{\times}X+k$. \end{proof}

Note that $G_{\infty}(f)=G(f)$ iff $f(X)$ is not of type a) or b) in 
theorem \ref{stich}.
 
{\noindent \bf Notations.}
In order to simplify the notations
for  $(1,b,1)\in S(f)$ we  denote by $\sigma_b$ the element
$\sigma_{1,b,1}$ and by $P_f(X,b)$ the corresponding polynomial $P_{1,b,1}(X)$. 
For $(1,y,1)$ and  $(1,z,1)$ in $S(f)$ we define the following function 
$\epsilon_f (y,z)=P_f(X,y)+P_f(X+y,z)-P_f(X,z)-P_f(X+z,y)$.

Now we can give the precise 
structure of the wild automorphism group:

\begin{corollary}\label{35}  Assume $g(C_f)>0$, then the element $\rho\in Z(G_{\infty,1}(f))$ and 
$G_{\infty,1}(f)=<\rho,\sigma_b>$. The commutation rule is given by
$[\sigma_y,\sigma_z]=\rho^{\epsilon_f(y,z)}$ where $\epsilon_f (y,z)=P_f(X,y)+P_f(X+y,z)-P_f(X,z)-P_f(X+z,y)\in \Z/p\Z$. Moreover we have the following
exact sequence $0\to <\rho>\simeq \Z/p\Z\to
G_{\infty,1}(f)\overset{\pi}{\rightarrow} 
k$ where $\pi (\sigma_y)=y\in k$ and the image of $\pi$ is 
finite dimensional as $\F_p$-vector space.
\end{corollary}

\begin{proof} Clear. \end{proof}

\section{Universal family - Modifications of covers -\\ Algorithm}

\subsection{Universal family}\label{41}

{\noindent \bf Notations.}
In order to be able to treat families of covers it is useful for a given
conductor $m$ prime to $p$ to work over the ring $A:=\F_p[t_i,1\leq i\leq m]$ and 
consider $f(X)=\sum_{1\leq i< m}t_iX^i+X^m$. 

A specialization homomorphism is an homomorphism  $\varphi :A\to k$,
where $k$ is an algebraically closed field of characteristic $p>0$; then 
$\varphi(f)(X)=\sum_{1\leq i< m}\varphi(t_i)X^i+X^m\in k[X]$.

Let $i<m$ and $(i,p)=1$ and denote by $n(i):=\max \{ n\in \N \ |\ ip^n<m\}$. The 
following lemma measures the defect for the polynomial $\varphi(f)$ to
be additive. For this we introduce $\Delta (f)(X,Y):=f(X+Y)-f(X)-f(Y)$.

\begin{lemma} \label{42} With the notations above there is a unique polynomial 
$F(X,Y)\in \bigoplus_{1\leq i< m,\ (i,p)=1}A[Y]X^{ip^{n(i)}}$ and a unique 
polynomial $P_f(X,Y)\in XA[Y][X]$ such that 
\begin{equation}\label{e1}
\Delta (f)(X,Y)=F(X,Y)+(Id-F)P_f(X,Y)
\end{equation}
The polynomial $P_f(X,Y)$ is characterized by the following:
\begin{equation} \label{e2}
P_f(X,Y)=(Id+F+...+F^n) \Delta (f) \mod  (X^{[\frac{m-1}{p}]})
\end{equation}
 for any 
$n$ such that $p^n> [\frac{m-1}{p}]$.
\end{lemma}

\begin{proof} Existence. Note that $\deg_X\Delta (f)=m-1$ and $\Delta f\in (X,Y)A[X,Y]$.
Let $(i,p)=1$ with $1\leq i <m$, and $0\leq j<n(i)$. For a monomial 
$a_{ip^j}(Y)X^{ip^j}$ of $\Delta (f)$ where $ a_{ip^j}(Y)\in A[Y]$ and total degree $<m$ we write 
$a_{ip^j}(Y)X^{ip^j}=(a_{ip^j}(Y))^{p^{n(i)-j}}X^{ip^{n(i)}}+(Id-F)(P_{ip^j}(X))$
with
$P_{ip^j}(X)=(Id+F+....+F^{n(i)-j-1})(a_{ip^j}(Y)X^{ip^j})$.

For the unicity we remark that if $P_f(X,Y)$ satisfies the formula \eqref{e1} then 
$\deg_X P_f(X,Y)\leq [\frac{m-1}{p}]$, so it is sufficient to prove the
formula \eqref{e2} in the lemma. We have the identity:

$(Id+F+...+F^{n-1}) \Delta (f)=(Id+F+...+F^{n-1}) F(X,Y)+(Id-F^n)P_f(X,Y)$

for any $n$. Now 
$F(X,Y)\in (X^{ [\frac{m-1}{p}]+1})$ as $ip^{n(i)}<m<ip^{n(i)+1}$ and
as $P_f(X,Y)\in XA[Y][X]$. Then for $p^n>[\frac{m-1}{p}]$ we obtain the
formula \eqref{e2}.
\end{proof}

\begin{defi} \rm Let $\varphi: A\to k$ be a specialization homomorphism. We denote by
the same letter the induced homomorphism on polynomials via the action
on the 
coefficients. Let $Ad_{\varphi(f)}(Y)\in k[Y]$ be the monic generator of the ideal of 
coefficients of  $\varphi (F)(X,Y)$ in $k[Y]$. As usual we denote by
$Z(Ad_{\varphi (f)}(Y))$ the set of zeroes in the algebraically closed
field $k$.
\end{defi} 

\begin{lemma} \label{44} Write $F(X,Y)=\sum_j a_j(Y)X^j=\sum_{1\leq i< m,\ (i,p)=1}a_{ip^{n(i)}}(Y)X^{ip^{n(i)}}$. 

i) The coefficient $a_{m-1}(Y)\in mY+(YA[Y])^p$, and $\deg a_{m-1}(Y)\leq (m-1)^2.$

ii) If $m=1+\el p^s$ where $\el >1,\ (\el,p)=1$, $s=v_p(m-1) > 0$,  
let $j_0=1+(\el -1)p^s$. Then 
$ j_0$ is prime to $p$ (also if $s=0$) and $n(j_0)=0$. Moreover the coefficient
$a_{j_0}(Y)=\el Y^{p^s}+...+2t_{j_0+1}Y$ if $p>2$ and if $p=2$ 
one has $a_{j_0}(Y)=\el Y^{p^s}+...+t_{j_0+2}Y^2$. 
If $(p,m-1)=1$ (so $s=0$ and $p>2$) one has  $a_{j_0}(Y)=(\el+1) Y=mY$.

\end{lemma}

\begin{proof}: For i) we remark that $m-1=\el p^s=\el p^{n(\el)}$ is the highest 
representative $<m$  of $\el$ modulo multiplication by a power of $p$.
As  $\Delta (X^m)=mYX^{m-1}+\mbox {\rm lower degree terms}$ and as lower degree monomials give contributions in $YA[Y]$ that one needs to raise to some $p$-power, the result follows. Concerning the degree, we remark that $\Delta(f)(X,Y)=
\sum_{i\leq m-1} \delta_i(Y)X^i$ and $\deg \delta_i(Y)\leq m-1$. Now write $i=jp^{n(j)}$ then the contribution of $\delta_i(Y)X^i$ in $F(X,Y)$ is $(\delta_i(Y)X^i)^{p^{n(j)}}$ 
whose $Y$-degree is $\leq (m-1)^2$.

For ii) we remark that $j_0>\frac{m-1}{p}+1$ and $j_0$ is prime to
$p$.  So the coefficient of the monomial $X^{j_0}$ in $F(X,Y)$ is the same as that of $\Delta(f)$ and 
so equal $\sum_{j_0<i\leq m,(i,p)=1}\binom{i}{j_0}t_iY^{i-j_0}=
\el Y^{p^s}+ \dots +(j_0+1)t_{j_0+1}Y$. 
The result follows. 
\end{proof}

\begin{proposition}\label{additive}  Let $\varphi: A\to k$ be a
specialization homomorphism. Then 
$\Ad_{\varphi (f)}(Y)$ is a separable and  additive polynomial and
$Z(\Ad_{\varphi (f)})=\{y\in k\ |\ \Delta (\varphi\circ f)(X,y) \in (Id -F)k
[X]\}$. Moreover if $y\in Z(\Ad_{\varphi (f)})$ then there is a unique 
$P(X)\in Xk[X]$ such that $\Delta (\varphi( f))(X,y)=(Id -F)P(X)$ and 
$P(X)=\varphi(P_f)(X,y)$. Then $Z(\Ad_{\varphi (f)}) =\{y\in k\ |\ \sigma_y\in G_{\infty,1}(f)\}$; in 
particular $| G_{\infty,1}(\varphi (f))|=\deg \Ad_{\varphi(f)}(Y)$ and for 
$y,z\in Z(\Ad_{\varphi (f)})$ the commutation rule for $\sigma_y,\sigma_z\in
G_{\infty,1}(\varphi (f))$ is 
$\epsilon_{\varphi(f)} (y,z)=\varphi(P_f)(X,y)+\varphi(P_f)(X+y,z)-\varphi(P_f)(X,z)-\varphi(P_f)(X+z,y)$.
\end{proposition}

\begin{proof} $\Ad_{\varphi(f)}(Y)$ is separable because we know from lemma
\ref{44} i) 
it  divides the polynomial $a_{m-1}(Y)\in mY+(YA[Y])^p$ which is separable.
Now to prove it is additive, it suffices to show its set of  roots is
stable under addition. We first remark that $Z(\Ad_{\varphi (f)})
\subset \{y\in k\ |\ \sigma_y\in G_{\infty,1}(f)\}$. For the reverse inclusion
we remark that in the equality
\begin{equation} \label{e3}
\Delta \varphi(f)(X,y)=(Id -F)P(X)
\end{equation}
where $P(X)\in k[X]$, we can assume that 
$P(X)\in Xk[X]$ as $\Delta (\varphi(f))(0,y)=0$. 
Note that $\deg P\leq [\frac{m-1}{p}]$; then for $n>>0$
$(Id+....+F^n)\Delta (\varphi(f))(X,y)=P(X) \mod X^{[\frac{m-1}{p}]+1}$,
so $P(X)=\varphi(P_f)(X,y)$. Then from \ref{42} and \eqref{e3} it follows that $\varphi(F)(X,y)=0$, i.e.\ $\Ad_{\varphi (f)}(y)=0$.

Let $y,z\in k$, such that $\Delta (\varphi(f))(X,y), \Delta (\varphi(f))(X,z)\in (Id -F)k[X]$. We have the following general identity
$\Delta (\varphi(f))(X,y+z)=\Delta (\varphi(f))(X+z,y)+\Delta
(\varphi(f))(X,z)-\Delta (\varphi(f))(y,z)$ and for our choice of
$y,z$ each term on the right hand side 
is in $(Id-F)k[X]$, so $\Ad_{\varphi (f)}(y+z)=0$. 
\end{proof}

\begin{corollary} We denote by $\Id: A\to \Fr A$ the inclusion
homomorphism. Then 

A. i) If $p=2$ and $m=3$ then $\Ad_{\Id (f)}(Y)=Y^4+Y$.

ii) If $p=2$ and $m=5$ then $\Ad_{\Id (f)}(Y)=Y^{16}+t_3^4Y^8+t_3Y^2+Y$.

iii) If $p=3$ and $m=4$ then $\Ad_{\Id (f)}(Y)=Y^9+2t_2^3Y^3+Y$.

B. Outside case A. One has $\Ad_{\Id (f)}(Y)=Y$.
Moreover $\exists D(\underline t)\in \F_p[\underline
t]-\{0\}$ such that for any specialization $\varphi: A\to k$ with
$\varphi ( D(\underline t))\neq 0$ one has $\Ad_{\varphi}(Y)=Y$ and so
$G_{\infty,1}\simeq \Z/p\Z$
\end{corollary}

\begin{proof}

A. A direct calculation gives the formulas.

B.  If $p$ doesn't divide $m-1$ the result follows from Lemma \ref{44} ii).
Therefore we can write $m=1+\el p^s$ with $(\el,p)=1,\ s>0$. We distinguish two
cases.

1. $\el >1$; let $j_0=1+(\el -1)p^s$, write $F(X,Y)=\sum_ja_j(Y)X^j$ then $a_{j_0}(Y)=\el Y^{p^s}+ 
\mbox{\rm lower degree terms }$.
The ring $A$ is factorial and $F(X,Y)\in A[Y][X]$, then the content of 
$F(X,Y)$ in $A[Y][X]$ is unitary and so equal to $\Ad_{\Id}(Y)$. This 
show that $\Ad_{\Id}(Y)\in A[Y][X]$ is a unitary polynomial. Let
$\varphi: A\to k$ be a specialization morphism then $\varphi (\Ad_{\Id})$ 
is unitary and $\deg \varphi (\Ad_{\Id})=\deg \Ad_{\Id}(Y)$. Now we 
remark that $\varphi (\Ad_{\Id})$
 divides $\varphi (a_{m-1}(Y))\in mY+(Yk [Y])^p$ and so is still 
separable in $k[X]$ (cf. Lemma \ref{44} i).
Let us now consider the specialization homomorphism  $\varphi_0:A\to 
\F_p$ defined by $\varphi_0(t_i)=0$ for
$i=1,...,m-1$. Then $\varphi_0(f)(X)=X^m$ and $\Ad_{\varphi_0}(Y)$ divides $ \el 
Y^{p^s}$
(see the proof of Lemma \ref{44} ii). It follows that  $\Ad_{\Id}(Y)=Y$.

2. $\el =1$, then $m=1+p^s$; with $s>0$. We show that $\Ad_{\Id}(Y)$ is
unitary: A look at the monomials in $\Delta (f)(X,Y)$ shows that the 
highest
contribution in $a_{m-1}(Y)$ comes from $(X+Y)^{p^s}$ and more precisely 
from
the linear term $Y^{p^s}X$ which we need to raise to the power $p^s$. 
So finally $a_{m-1}(Y)=Y^{p^{2s}}$+ lower degree monomials. Now like in 
case 1 we look for good specialization morphisms.
For $p>3$, we consider $\varphi_0(f)(X)=X^3+X^m$, then $\Ad_{\varphi_0}(Y)|
\varphi_0(a_{2p^{s-1}}(Y))=Y^{p^{s-1}}$, we conclude as in case 1.

For $p=2$ and $s>2$, we consider $\varphi_0(f)(X)=X^7+X^m$, then 
$\Ad_{\varphi_0}(Y)|
\varphi_0(a_{6p^{s-2}}(Y))=Y^{p^{s-2}}$, we conclude as in case 1.

For $p=3$ and $s>1$, we consider $\varphi_0(f)(X)=X^5+X^m$, then 
$\Ad_{\varphi_0}(Y)|
\varphi_0(a_{4p^{s-1}}(Y))=(2Y)^{p^{s-1}}$, we conclude as in case 1.

In order to exhibit a polynomial $D(\underline t)$ we recall  that $Y$ is 
the content of the polynomial $F(X,Y)=\sum_ja_j(Y)X^j$
 in $A[Y][X]$; so there are $b_j(Y)\in \Fr A[Y]$ such that
$\sum_jb_j(Y)a_j(Y)=Y$.

Any  $D(\underline t)\in \F_p[\underline t]-\{0\}$ such that 
$D(\underline t)
b_j(Y)\in A[Y]$ for all $j$ works.

\end{proof}

\begin{remark} \rm The method used here is a special case of
\cite[Theorem 1.11 and  lemma 1.12]{De-Mu}. As $A$ is a UFD it is a
natural question to ask for the best $D(\underline{t})$.

\end{remark}

\subsection{Modifications of covers} 

In this paragraph $k$ is an algebraically closed field ($\car p>0$).
In the previous paragraph we fixed the conductor which is equivalent
to fixing
the genus of the cover. We study here how for $f(X)\in k[X]$ the additive polynomial $\Ad_f(Y)$ changes through
natural algebraic transformations which don't necessarily preserve the genus.

\begin{defi} \rm
Let $S(X)$ be an additive polynomial $\in k[X]$. We say that the cover 
$C_{f\circ S}$ is a {\it modification of type 1} of $C_f$. Although $f(X)$ is 
reduced, in general $f\circ S$ is not reduced. 

\end{defi}

\begin{proposition} \label{49} Let $S(X)\in k[X]$ be a separable and
additive polynomial. 
Then $\Ad_{f}(S(Y))$ divides  $\Ad_{\red(f\circ S)}(Y)$. 
Further for $y,z\in Z(\Ad_{f}(S(Y)))$, we have $y,z\in Z(\Ad_{\red(f \circ S)}(Y))$ and 
$\epsilon_f (S(y),S(z))=\epsilon_{\red (f\circ S)} (y,z)$ (see \ref{additive}).
\end{proposition}

\begin{proof}
One can write $\Delta (f)(X,Y)=\Ad_f(Y) G(X,Y)+P_f(X,Y)-P_f(X,Y)^p$, where 
$G(X,Y)\in k[Y][X]$ has content equal to $1$, and $P(X,Y)\in k[Y][X]$.
Then $\Delta (f\circ S)(X,Y)=\Delta (f)(S(X),S(Y))=$
$$\Ad_f(S(Y)) G(S(X),S(Y))+P_f(S(X),S(Y))-P_f(S(X),S(Y))^p.$$ Let $y\in
Z(\Ad_f(S(Y)))$, 
then $\Delta (f\circ S)(X,y)=P(S(X),S(y))-P(S(X),S(y))^p$ and as
$\red(f\circ S)(X)=f\circ S(X)+Q(X)-Q(X)^p$ for some $Q(X)\in k[X]$ we can write$\Delta (\red{(f\circ S)})(X,y)=(\Id-F)(P(S(X),S(y))+\Delta(Q)(X,y))$ which by proposition \ref{additive} gives $\Ad_{f\circ S}(y)=0$. The 
divisibility follows as $S(Y)$ and so $\Ad_{f}(S(Y))$ is separable. Moreover 
by the same proposition for $y\in Z(\Ad_{f}(S(Y)))$, we have $P_{\red(f\circ S)}(X,y)=P(S(X),S(y))+\Delta(Q)(X,y)$. We remark that $\Delta(Q) (X,y)+\Delta(Q) (X+y,z)-\Delta(Q) (X,z)-\Delta(Q) (X+z,y)=0$; it follows that for $y,z\in Z(\Ad_{f}(S(Y)))$ 
$\epsilon_f (S(y),S(z))=\epsilon_{\red (f\circ S)} (y,z)$.
\end{proof}

\begin{remark} \rm The divisibility can be strict. For example let $f(X)=X^{1+p^n}$,
then $\Ad_f(Y)=Y+Y^{p^{2n}}$. Let $S(X)=X+X^p$, then $f\circ S(X)=(X+X^p)(X^{p^n}+X^{p^{n+1}})=X^{1+p^n}+X^{1+p^{n+1}}+X^{p+p^n}+X^{p+p^{n+1}}$ and so $\red(f\circ S)=X^{1+p^{n-1}}+2X^{1+p^n}+X^{1+p^{n+1}}$.
Further, by Proposition \ref{415} below,
$\Ad_{f\circ S}(Y)=Y+...+Y^{p^{2(n+1)}}$ and $\Ad_{f}(S(Y))=Y+Y^p+Y^{p^{2n}}+Y^{p^{2n+1}}$. 
\end{remark}

\begin{defi} \rm 
Let $f(X),g(X)\in k[X]$, we assume that none of $\red(f),\red(g),\red(f+g)$ is 
zero. If $\Ad_{\red(g)(Y)}|\Ad_{\red(f)}(Y)$ we say that the cover $C_{f+g}$ is a 
{\it modification of type 2} of $C_f$.  

\end{defi} 

\begin{proposition} \label{412} Let $f,g$ as above and assume that $C_{f+g}$ be a modification of type 2 of
$C_f$. Then $\Ad_{\red(g)}(Y)|\Ad_{\red(f+g)}(Y)$.
Further for $y,z\in Z(\Ad_{\red(g)}(Y))$ we have $y,z\in Z(\Ad_{\red(f+g)}(Y))$ and 
$\epsilon_{\red(f+g)} (y,z)=\epsilon_{\red(f)} (y,z)+\epsilon_{\red(g)} (y,z)$ (see \ref{additive}).
\end{proposition}

\begin{proof} Clear from Proposition \ref{additive}. \end{proof}

\begin{proposition} \label{413} Let $k$ be an algebraically closed field of $\car p>0$ and $f(X)\in k[X]$ whose degree $\deg f:=m$ is prime to $p$. Let $g(f)=\frac{(m-1)(p-1)}{2}$, the genus of $C_f$; we assume $g(f)>1$. We write 
$m=1+\el p^s$ with $\el \geq 1,\ (\el,p)=1 $ and $s > 0$. 

i) Let $\el >1$ and $p=2$, then  $|G_{\infty,1}(f)|\leq 2^s$ and the
ratio $\frac{|G_{\infty,1}(f)|}{g(f)}\leq \frac{2}{3}$. The two
inequalities  are equalities  for 
$f(X)=\red {(( X+X^{2^{s-1}})^7)}$ and $s>1$.

ii) Let $\el >1$ and $p>2$, then $|G_{\infty,1}(f)|\leq p^{s+1}$ and
the ratio $\frac{|G_{\infty,1}(f)|}{g(f)}\leq \frac{2}{\el}\frac{p}{p-1}\leq \frac{p}{p-1}$. We 
have equality for $f(X)=X^{1+2p^s}-X^{2+p^s}$.

iii) Let $\el =1$, i.e.\ $m=1+p^s$, then $|G_{\infty,1}(f)|\leq p^{2s+1}$,
i.e.\ ratio $\frac{|G_{\infty,1}(f)|}{g(f)}\leq 2p^s\frac{p}{p-1}$ and we 
have equality for $f(X)=X^{1+p^s}$.

\end{proposition}

\begin{proof}
i) and ii). We work with the universal family defined by $f\in A[X]$ of degree $m$ (cf.\ref{41}). 
Let $j_0=1+(\el -1)p^s$. We have seen in Lemma \ref{44} that 
$a_{j_0}(Y)=\el Y^{p^s}+...+2t_{j_0+1}Y$. As $\Ad_{\varphi}(Y)$ is a separable 
and additive polynomial by Proposition \ref{additive}, it follows that 
$\deg \Ad_{\varphi}(Y)\leq \deg a_{j_0}(Y)=p^s$ if $p>2$ and  $\deg \Ad_{\varphi}(Y)\leq(1/2) \deg a_{j_0}(Y)=2^{s-1}$ if $p=2$.

For the equality of the bound we give examples. 

Let $p=2$, $(\el,2)=1$ and $f(X)=X^{1+\el 2}$. 
In this case it is an easy consequence of Lemma \ref{44} that $\Ad_f(Y)=Y$.
Let $S(X)=X+X^{2^{s-1}}$ where $s>1$. Then $f\circ S(X)=(
X+X^{2^{s-1}})(X^2+X^{2^s})^\el=X(X^2+X^{2^s})^\el+X^{2^{s-1}}(X+X^{2^{s-1}})^\el\mod
(\Id -F)k[X]$ has conductor $1+\el 2^s=m$. From Proposition \ref{49} we know
that $\Ad_f(S(Y))=Y+Y^{2^{s-1}}|\Ad_{ f\circ S}(Y)$ and from the first part of
the Proposition $\deg \Ad_{ f\circ S}(Y)\leq 2^{s-1}$. So we have equality
$\Ad_{ f\circ S}(Y)=Y+Y^{2^s}=\Ad_f(S(Y))$. Then the ratio
$\frac{|G_{\infty,1}(f\circ S)|}{g(f)}= \frac{2}{\el}\leq
\frac{2}{3}$. So we have equality for $l=3$, i.e.\ $f(S(X))=(X+X^{2^{s-1}})^7$.

Let $p>2$ we give an example for $\el=2$. 
As $j_0+1=2+(\el-1)p^s=2+p^s$ we need $t_{j_0+1}\neq 0$. So we consider 
$f(X)=X^{1+2p^s}-X^{2+p^s}$ and show that $\Ad_f(Y)=Y^{p^s}-Y$. As  
from the first part of the Proposition $\deg \Ad_{f}(Y)\leq p^s$ it suffices to
prove a divisibility. Let $y\in k\ |\ y^{p^s}=y$, then 
$\Delta (f)(X,y)=(X+y)(X^{p^s}+y)^2-(X+y)^2(X^{p^s}+y)-f(X)-f(y)=
(X+y)(X^{2p^s}+2yX^{p^s}+y^2)-(X^2+2yX+y^2)(X^{p^s}+y) -f(X)-f(y)=0\mod (\Id-F)k[X]$

iii) see Lemma \ref{44} i). 
\end{proof} 

\begin{remark} \rm In general for a given curve $C$ in order to bound the automorphism group one considers $\el>2 $ prime  $\neq p$ and uses the into homomorphism 
$$\Aut C \hookrightarrow \Gl(2g(C),\F_{\el})$$ where $g(C)$ is the genus. Let us 
for example consider the case $p=2$ and $m=1+2^k$. Then $2g(C_f)=2^k$, and 
so $|G_{\infty,1}(f)|$ divides the cardinality of a $2$-Sylow of $ \Gl(2^k,\F_{\el})$ which is equal to $2^{g(C)F(1)-1}$,
where $F(1)=1+v_2(\el^2-1)+v_2(\el -1)$ (see \cite {C-F} for the structure 
of $2$-Sylow subgroups of $\Gl(2^k,\F_{\el})$. From this one also
obtains bounds for $\Aut C$ but they are far from being as good as
the ones in Proposition \ref{413}.
\end{remark}

\subsection{Extraspecial groups} \label{extra}
We recall some basic facts concerning extraspecial groups. We refer to 
\cite{Hu} and \cite{Su} for the structure of finite $p$-groups.

We saw in Corollary \ref{35} that the groups $G_{\infty,1}(f)$ belong to  the class $C_1$ of $p$-groups:

$$C_1:=\{ G\ |\ \exists C\  p-{\mbox{\rm cyclic}}\subset Z(G)\ 
{\mbox{\rm such that}}\
G/C\ \mbox{\rm is } p \mbox{\rm -elementary abelian}\}$$

Those non abelian groups in $C_1$ for which $Z(G)$ is itself $p$-cyclic are 
called extraspecial. In particular a non abelian group $G$ of order $p^3$ is 
 extraspecial:

If $p=2$, $G$ is isomorphic to the dihedral group $D_8$ or the 
quaternion group $Q_8$. These two groups have exponent $2^2$.
If $p>2$, $G$ is isomorphic to one of the two groups:

$E(p^3)=<x,y|x^p=y^p=[x,y]^p=1,[x,y]\in Z(E(p^3))>$ with exponent $p$, or

$M(p^3)=<x,y|x^{p^2}=y^p=1,y^{-1}xy=x^{1+p}>$ with exponent $p^2$.

More generally let $G$ be an extraspecial $p$-group then (see \cite{Su} Th. 4.18) $|G|=p^{2n+1}$ for some $n>0$ and the following 4 types occur.

I. If $\exp G=p$, then $p>2$ and $G$ is a central product of $n$ groups $E(p^3)$.

II. If $\exp G=p^2$ and $p>2$ then $G$ is a central product of $M(p^3)$ and $n-1$ groups   $E(p^3)$.

III. If $p=2$, then either 

a. $G$ is a central product of $n$ groups $D_8$. Let $q_a(n)$ (resp.
$d_a(n)$) the number of elements with order $4$ (resp. $\leq 2)$.

Or b. $G$ is a central product of a group $Q_8$ and $n-1$ groups  $D_8$. 
Let $q_b(n)$ (resp. $d_b(n)$) the number of elements with order $4$ (resp. $\leq 2)$. 

We have the following:

$q_a(n)=-2^n+2^{2n}$ and so $d_a(n)=2^n+2^{2n}$

$q_b(n)=2^n+2^{2n}$ and so $d_b(n)=-2^n+2^{2n}$

In each case, the structure of the central product is uniquely determined by the structure of the factors. In the sequel when speaking of the isomorphism type of an extraspecial group we will refer to the  4 types above, namely  I, II, III.a and III.b.

\begin{proposition} \label{415} Let $n>1$, $A_n:= \F_p^{alg}[t_i,0\leq i\leq n-1]$, and $f:=t_0X^{1+p^0}+t_1X^{1+p}+...+t_{n-1}X^{1+p^{n-1}}+X^{1+p^n}\in A_n[X]$. Let $k$ be an 
algebraically closed field and $\varphi: A_n\to k$, a specialization
homomorphism. 
Then $\Ad_f(Y)=\sum_{0\leq i\leq n}(t_i^{p^{n-i}}Y^{p^{n-i}}+
t_i^{p^n}Y^{p^{n+i}})=t_nY+...+t_n^{p^n}Y^{p^{2n}}$ and
$|G_{\infty,1}(\varphi(f)) |=p^{2n+1}$. Further 

i) If $p>2$, then $G_{\infty,1}(\varphi(f))$ is an extraspecial group of type I.

ii) If $p=2$, and $t_0=0$, then $G_{\infty,1}(\varphi(f))$ is an extraspecial group of type III.b.

\end{proposition}

\begin{proof}
Write $f_i= t_iX^{1+p^i}$, with $t_n=1$ and $t_0=0$ for $p=2$. Then $\Delta (f)(X,Y)=\sum_{0\leq i\leq n}\Delta (f_i)$, where $\Delta (f_i)= A_i+B_i$, $A_i=t_iYX^{p^i}$, $B_i=t_iY^{p^i}X$. We can write $\Delta (f_i)=(t_i^{p^{n-i}}Y^{p^{n-i}}+
t_i^{p^n}Y^{p^{n+i}})X^{p^n}+P_i-P_i^p$. $P_i=A_i+...+A_i^{p^{n-i-1}}+B_i+...+
B_i^{p^{n-1}}$. 
One obtains $\Ad_f(Y)= \sum_{0\leq i\leq n}(t_i^{p^{n-i}}Y^{p^{n-i}}+
t_i^{p^n}Y^{p^{n+i}})=t_nY+...+t_n^{p^n}Y^{p^{2n}}$. The polynomial 
$P(X,Y)=\sum_{0\leq i\leq n}P_i(X,Y)$ is here an additive polynomial, so
if $\sigma_{c,Y}(W)=W+P(X,Y)+c$, for $c\in \F_p$, then 
$\sigma_{c,Y}^p(W)=W+P(X,Y)+P(X+Y,Y)+...+P(X+(p-1)Y,Y)=W+\frac {p(p-1)}{2}
P(Y,Y)$ and so if $p>2$, $\sigma_{c,Y}^p=\Id $ and the exponent  of
$G_{\infty,1}(\varphi(f))$ is $p$. Now  if $p=2$, $\sigma_{c,Y}^2(W)=W+P(Y,Y)$ 
and the exponent  of
$G_{\infty,1}(\varphi(f))$ is $2^2$.

Let $y,z\in Z(\Ad_{\varphi (f)}(Y))$, then $\epsilon_f(y,z)=P(X,z)+P(X+z,y)-
P(X,y)-P(X+y,z)=P(z,y)-P(y,z)=\sum_{0\leq i\leq n}\sum_{n-i\leq j\leq n-1}
(B_i(z,y)-A_i(z,y))^{p^j}=...-t_n^{p^{n-1}}y^{p^{n-1}}z^{p^{2n-1}}$. If 
$y\neq 0$, is such that $\sigma_y$ is in the center this polynomial in $z$ 
should have $p^{2n}$ 
roots ; so the center is $<\rho>$ where $\rho(X)=X$ and $\rho(W)=W+1$.
We conclude that  $G_{\infty,1}(\varphi(f))$ is an extraspecial group of 
order $p^{2n+1}$ and exponent $p$ for $p>2$ and exponent $p^2$ for $p=2$.

In the case $p=2$ this is not yet sufficient to determine the isomorphism class
of this group; we need some evaluation of the number of elements of order $\leq 2$.
Let us first consider the case $f_0(X)=X^{1+2^n}$. We have the following 
parametrization for the elements in $G_{\infty,1}(f_0)$.

For $y\in Z(\Ad_{f_0}(Y))=Z(Y+Y^{2^{2n}})$ and $c\in Z/pZ$, $\sigma_{y,c}(X)=X+y$ and $\sigma_{y,c}(W)=W+P(X,y)+c$, where $P(X,Y)=Y^{2^n}X+Y^{2^{n+1}}X^2+...+Y^{2^{2n-1}}X^{2^{n-1}}$. 
Let $d(n)$ be the number of elements of order $\leq 2$; 
they correspond to those $y$ such  that $P(y,y)=0$ (this condition doesn't depend
on $c$). One has $P(Y,Y)=S+S^2+...+S^{2^{n-1}}$, where $S=Y^{2^n+1}$ and so $Q(Y):=P(Y,Y)Y^{-2^n}$ is a separable polynomial of degree $2^{2n-1}+2^{n-1}-2^n$.
Moreover $P(Y,Y)^2-P(Y,Y)=S^{2^n}-S=Y^{2^n}(Y^{2^{2n}}+Y)$. So the roots of $Q(Y)$ are simple and among those of $Y^{2^{2n}}+Y$. It follows that $d(n)=2(2^{2n-1}+2^{n-1}-2^n)=2^{2n}-2^n=d_b(n)$.

Now let us consider the general case. The result follows from the previous special case if we remark that the number of elements of order 4 can only decrease 
via a specialization $\varphi: A_n\to k$. And as $q_b(n)>q_a(n)$, only extraspecial groups of type b) can degenerate in type b). 

\end{proof}

\begin{remark} \rm 
It was surprising for us when we saw that a special case appears in literature in connection with coding theory; namely \cite{vdG-vdV} consider  the specialization morphisms  $\varphi$ taking values in $\F_q$ where $q=p^n$. The additive polynomial $\Ad_{\varphi (f)}$ is
then their polynomial $E_n(Y)=R(Y)^{p^n}+\sum_{0\leq i\leq n}(t_iY)^{p^{n-i}}$ where $R(X)=\sum_{0\leq i\leq n}t_iY^{p^i}$. The zeroes are interpreted as the $\F_p$ vector space which is the kernel of the $\F_p$-bilinear form 
$\Tr_{F_{q}/\F_p}(xR(y)+yR(x))$ (which is symmetric 
if $p>2$ and alternating if $p=2$).  
In case $p=2$ they prove a  factorization of  $E_n(Y)=YE_n^{-}(Y)E_n^{+}(Y)$
which corresponds to a partition of the roots depending on the order
of the corresponding automorphism of $C_{f}$. Such a decomposition works in general.  
\end{remark}

\subsection{Application to the moduli space of curves} 

We keep the notations of Proposition \ref{415}.
For a fixed $n>1$ let $g_n=\frac{p^n(p-1)}{2}$, $A_n:=
\F_p^{alg}[t_i,0\leq i\leq n-1]$
 and
$f:=t_0X^{1+p^0}+t_1X^{1+p}+...+t_{n-1}X^{1+p^{n-1}}+X^{1+p^n}\in
A_n[X]$ 
with $t_0=0$ for $p=2$.
Let $\theta \in \F_p^{alg} $ be a primitive $(p-1)(p^n+1)$-th root of $1$.
Then $\Theta: (t_0,...,t_{n-1})\to (\theta ^{p^0-p^n}t_0,\theta ^{p^1-p^n}t_1,...,\theta ^{p^{n-1}-p^n}t_{n-1})$ induces  an $\F_p^{alg}$-automorphism of $\Spec
A_n$ of order $p^n+1$. Let $B_n:=A_n^{<\Theta>}$ be the quotient of the affine space by the cyclic group of automorphisms $<\Theta>$.   
The structural morphism $\pi: C_{f}\to \Spec   A_n$ is a family of
curves of genus $g_n$. 
Moreover by Proposition \ref{33}
two specialization morphisms $\varphi_i: A_n\to k$, for $i=1,2$
will give  isomorphic $k$-curves iff $\exists c\in \F_p^{\times}$, 
$(a,b)\in (k^{\times},k)$ such that 
$\sum_{0\leq i\leq n-1}\varphi_2(t_i)(aX+b)^{1+p^i}+(aX+b)^{1+p^n}
=c(\sum_{0\leq i\leq n-1}\varphi_1(t_i)X^{1+p^i}+X^{1+p^n}) \mod (\Id-F)k[X]$,
i.e.\ $a^{1+p^n}=c\in \F_p^{\times}$ and for $i\geq 0$ one has $\varphi_2(t_i)=
ca^{-(1+p^i)}\varphi_1(t_i)=a^{p^n-p^i}\varphi_1(t_i)$ (note that for $p=2$
we assumed $t_0=0$ and so $f$ is reduced). Finally this shows that the two specialization
morphisms are in the same orbit under the action of the group $\Theta$. 
 By the definition of the coarse moduli space $M_{g_n}$ we deduce from
the existence of the family $ C_{f}\to \Spec A_n$  a    
map from  $\Spec A_n$ to $ M_{g_n}$ which factorizes through $\Spec
B_n$ in an injective morphism. The image is an algebraic subset of $
M_{g_n}$ of dimension that of  $\Spec B_n$. A measure of the size of
families of curves which are \'etale covers of the affine line and
with given extrasspecial group of type I and order $p^{2n+1}$ as an
automorphism group is given by the dimension of this  image which is
$O(\mbox{\rm log} (g_n))$. 

It is remarkable that these families for varying $n$ can be characterized by the following Hurwitz-type bound.

\begin{proposition} \label{417} Let $k$ be an algebraically closed
field of $\car p>0$ and $f(X)\in Xk[X]$ 
a polynomial of degree $m:=\deg f$ prime to $p$. We assume $f$ is reduced.
If $\frac{|G_{\infty,1}(f)|}{g(f)}>\frac{p}{p-1}$ 
($2/3$ for $p=2$ ) then $f(X)=XR(X)$, where $R(X)$ is an additive polynomial. Moreover if $\deg f=1+p^n$, then $\frac{|G_{\infty,1}(f)|}{g(f)}
=2p^n\frac{p}{p-1}$.
\end{proposition}

\begin{proof} We will show only the case $p>2$ and point out that for
$p=2$ a similar argument works.
The proof works by elimination of bad monomials. 
We saw in Proposition \ref{413} that only for 
$m=1+p^s$ with $s>0$ the ratio $\frac{|G_{\infty,1}(f)|}{g(f)}$ can be $>\frac{p}{p-1}$. Now we show that any other monomial in $f(X)$ has exponent 
$1+p^t$ with $t<s$. 

Let us assume that this is not the case and denote by $X^a$ the 
monomial of highest degree which is not of the above form. 
We first assume that $p$ doesn't divide $a-1$ and consider the integer
$k$ such that $p^{k-1}<a-1<p^k$. Then $k\leq s$ and $p^{s-1}<(a-1)p^{s-k}<p^s$. Then $(*)X^a$ 
has a contribution in $F(X,Y)$ which is equal to  $(*)a^{p^{s-k}}Y^{p^{s-k}}
X^{(a-1)p^{s-k}}+\mbox {\rm lower degree terms}$. Moreover as for $\alpha >0$,
$p^{\alpha}(a-1)>(a-1)$ it follows that $(*)a^{p^{s-k}}Y^{p^{s-k}}$ is exactly 
the coefficient of $X^{(a-1)p^{s-k}}$ in $F$. This would imply that
$|G_{\infty,1}(f)|\leq p$ and the ratio
$\frac{|G_{\infty,1}(f)|}{g(f)}\leq \frac{2p}{(p-1)p^s}\leq
\frac{p}{(p-1)}$ - a contradiction!

Let us now assume that $a-1=\el p^t$, where $\el >1$ and $(\el, p)=1$. Let 
$j_0=1+(\el -1)p^t$ and say $p^{k-1}<j_0<p^k$. Then
$p^{s-1}<j_0p^{s-k}<p^s$ and
the monomial $X^a$ contributes to $F(X,Y)$ in the monomial 
$X^{j_0p^{s-k}}$ the term $(Y^{p^t}+....+t_{j_0+1}Y)^{p^{s-k}}X^{j_0p^{s-k}}$.
Note that $pj_0<a$ iff $\el<\frac{p}{p-1}-\frac{1}{p^t}$ which is not the case;
so any other contribution in $F$ in the monomial $X^{j_0p^{s-k}}$ can only
occur from a monomial $X^b$ with $j_0<b<a$. Now such a contribution will 
be $(\binom{b}{j_0}Y^{b-j_0})^{p^{s-k}}$ whose degree is $(b-j_0)p^{s-k}
<(a-j_0)p^{s-k}=p^{t+s-k}$, so finally we get that 
$\deg Ad_f(Y)\leq p^{t+s-k}$.
From this we get the ratio  $\frac{|G_{\infty,1}(f)|}{g(f)}\leq
\frac{2p^{t+s-k+1}}{p^s(p-1)}$ and so $1<2p^{t-k}$, i.e.\ $t\geq k$. On
the other hand $(\el-1)p^t<p^k$, so
$(\el -1)<\frac{1}{p^{t-k}}$ - a contradiction.
\end{proof}

\subsection{Realization of the other extraspecial groups} \label{4.5}

Case $p>2$.

Let $G$ be an extraspecial group of type II. We will use Witt vectors of length
2 and modifications of covers.

Let $c(X,Y)=\frac{(X+Y)^p-X^p-Y^p}{p}=\sum_{1\leq i\leq  p-1}\frac{(-1)^{i-1}}{i}X^iY^{p-i}$ and 
$f_1(X):=c(X^p,-X)=\sum_{1\leq i\leq  p-1}\frac{1}{i}X^{p+(p-1)i}$.
A straightforward calculation in $W_2(\F_p)$ shows that:

$\Delta (f_1)(X,Y)=c((F-\Id)X,(F-\Id)Y)+(F-\Id)c(X,Y)$. As $c((F-\Id)X,(F-\Id)Y)
=(Y^p-Y)(X^p-X)^{p-1}+\mbox {\rm lower degree terms}\in (Y^p-Y)\F_p[X,Y]$ it 
follows that $\Ad_{f_1}(Y)=Y^p-Y$ and for $y\in Z(\Ad_{f_1}(Y))$ (i.e.\ $y\in 
\F_p$) $P(X,y)=-c(X,y)=\sum_{1\leq i\leq  p-1}\frac{(-1)^i}{i}y^{p-i}X^i$.

Now we show that $G_{\infty,1}(f_1)$ is $p^2$-cyclic. 
Let  $\sigma_y(W)=W+P(X,y)$. Then $\sigma_y^p(W)=W-\frac{(-1)^{p-1}}{p-1}y^{p-(p-1)}=W+y$, i.e.\ $\sigma_y^p=\rho^y$, so $G_{\infty,1}(f_1)=<\sigma_y>$ is 
$p^2$-cyclic. Note that in this case the conductor is $m=1+\el p$ with 
$\el=p-1$; then the ratio is $\frac{|G_{\infty,1}(f_1)|}{g(f_1)}=\frac{2p}{(p-1)^2}<\frac{p}{p-1} $ as $p>2$.
In order to get an extraspecial group of exponent $p^2$ we use modifications 
of $f_1$. Let $q:=p^n$, $f_2(X)=X^{1+q}$, $\theta $ a primitive 
$q^2-1$-th root of unity and $A(X):=\theta X^q-\theta^qX$. Then 
$S(X):=A(X)+A(X)^q+...+A(X)^{q/p}$ is an additive polynomial and $S(X)^p-S(X)=
A(X)^q-A(X)=\theta^q(X+X^{q^2})$.
Set $f(X):=f_1(S(X))+f_2(X)$.

\begin{proposition} \label{418} $G_{\infty,1}(f)$ is an extraspecial group of type II
and the ratio $\frac{|G_{\infty,1}(f)|}{g(f)}=\frac{2p}{(p-1)^2}<\frac{p}{p-1} $.

\end{proposition} 

\begin{proof} We have
$\Delta (f_1)(S(X),S(Y))=c(S(X)^p-S(X),S(Y)^p-S(Y))+(F-\Id)c(S(X),S(Y))=
c(A(X)^q-A(X),A(Y)^q-A(Y))+(F-\Id)c(S(X),S(Y))$ 
and $\Delta (f_2)(X,Y)=Y^qX+YX^q$.

We claim that $\Ad_{f_1\circ S}(Y)=\theta^{-q}\Ad_{f_1}(S(Y))=Y+Y^{q^2}$ and 
for $y+y^{q^2}=0$ one has $P_{f_1\circ S}(X,y)=c(S(X),S(y))$. To this
end notice that
$f_1(S(X))=c(S(X)^p,-S(X))=\sum_{1\leq i\leq p-1}\frac{1}{i}S(X)^{p+(p-1)i}$
has conductor $1+(p-1)q^2$, hence $\deg \Ad_{f_1\circ S}\leq  q^2$. 
As by Proposition \ref{49}
$\Ad_{f_1}(S(Y))=S(Y)^p-S(Y)=A(Y)^q-A(Y)=\theta^q(X+X^{q^2})|\Ad_{f_1\circ S}(Y)
$ we get the equality.

Note that $f_1\circ S$ and $f_2$ have the same additive polynomial, so
due to the property of second type modifications, we obtain 
$Y+Y^{q^2}|\Ad_f(Y)$. As $f$ and $f_1\circ S$ have the same conductor we get 
$\deg \Ad_f(Y)\leq q^2$. Finally we conclude $\Ad_f(Y)=Y+Y^{q^2}$
showing the claim.

Next we claim that $G_{\infty,1}(f)$ has exponent $p^2$. 
Let $y$ be such that $\Ad_f(y)=\Ad_{f_1\circ S}(y)=\Ad_{f_2}(y)=0$. Then 
$P_f(X,y)=P_1(X,y)+P_2(X,y)$, where $P_1(X,Y)=c(S(X),S(Y))$ and 
$P_2(X,Y)=-(y^qX+(y^qX)^p+...+(y^qX)^{q/p})$.

Note that $c(S(X),S(Y))=\sum_{1\leq i\leq p-1}\frac{(-1)^{i-1}}{i}S(y)^iS(X)^{p-i}$ and as $S$ is an additive polynomial we set $Z:=S(y)S(X)^p-S(y)^pS(X)$.
Then $\sum_{1\leq i\leq p-1}P_1(X+iy,y)=\Tr_{k(S(X))/k(Z)}c(S(X),S(y))=S(y)$ where
$k=\F_p^{alg}$. Note that $f_2$ induces an extraspecial group of exponent
$p$ so $\sum_{1\leq i\leq p-1}P_2(X+iy,y)=0$. 

Finally if $y^{q^2}+y=0$ and  $\sigma_y(W)= W+P_f(X,y)$ then
$\sigma_y^p(W)=W+S(y)$.

Now we show that the center is $<\rho>$. Let $y,z\in Z(\Ad_f(Y))$; then
$\epsilon_j(y,z)=P_j(X,z)+P_j(X+z,y)-P_j(X,y)-P_j(X+y,z)$ for $j=1,2$.

We have seen that $S(X)^p-S(X)=\theta^q(X^{q^2}+X)$, so if $y^{q^2}+y=0$
then $S(y)^p-S(y)=0$. We also saw that $\Ad _{f_1}(Y)=Y^p-Y$ and 
$G_{\infty,1}(f_1)$ is cyclic, so $\epsilon_{f_1\circ S}(y,z)=\epsilon_{f_1}
(S(y),S(z))=0$. Finally 
$\epsilon (y,z)=\epsilon_2(y,z)=-(z^qy+z^{qp}y^p+...+z^{q^2/p}y^{q/p})+
(zy^q+z^py^{pq}+...+z^{q/p}y^{q^2/p})$. For $z\neq 0$ this is a polynomial in 
$y$ of degree $q^2/p$, so it has at most $q^2/p$ roots and hence $z\in Z(G_{\infty,1}(f))$ iff $z=0$. 
 
\end{proof}

\begin{remark}\label{A} \rm 
One can follow the same method in order to produce a family of curves with 
automorphism group an extraspecial group of exponent $p^2$. 

Say $q=p^n$, 
$f_2=(t_0)^qX^{1+p^0}+(t_1)^qX^{1+p}+...+(t_{n-1})^qX^{1+p^{n-1}}+X^{1+p^n}$.
Then 
$\Ad_{f_2}(Y)=Y^{p^{2n}}+Y+(t_{n-1}^qY+{t_{n-1}}^{p^{2n-1}}Y^{p{2(n-1)}})^p+...
+(t_1^qY+t_1^{p^{n+1}}Y^{p^2})^{p^{n-1}}$.

Now we look for a $\theta $ such that $\theta\Ad_{f_2}(Y)\in
(F-\Id)k[X]$. Write 
$\theta =\alpha^{p^{2n}}$; then the condition can be simplified by using the 
equivalence $A^p=A \mod  (F-\Id)k[X]$. We get 

$$(\alpha^{p^{2n}}+\alpha)+(\alpha^{p^{2n-1}}t_{n-1}^{p^n}+\alpha^pt_{n-1}^{p})+...+(\alpha^{p^{n+1}}t_1^{p^n}+\alpha^{p^{n-1}}t_1^{p^{n-1}})=0.$$

Choose for $\alpha $ any root. We then can write $\theta \Ad_{f_2}(Y)=
S(Y)^p-S(Y)$ where $S(Y)$ is an additive polynomial, and this fact is essential for further calculations. Now $f=f_1\circ S+f_2$ gives a family of curves with 
automorphism group an extraspecial group of exponent $p^2$. 

\end{remark}

If $p=2$, it is not so easy to distinguish between the two classes of 
extraspecial groups. In fact we have realized those extraspecial groups which 
are a central product of a group $Q_8$ and $n-1$ groups  $D_8$, i.e.\ the so called
type III.b. A realization of type III.a. will be a consequence of the following
general method.

\subsection{Saturated subgroups of extraspecial groups and their realization}

Now we give the main result of this paper which describes the set of groups $G_{\infty,1}(f)$. 
Let us first define 2 sets of isomorphism classes of $p$-groups.

$$C_1:=\{ G\ |\ \exists N\  p{\mbox{\rm -cyclic and normal subgroup}}\subset G\ |\ 
G/N\ \mbox{\rm is } p \mbox{\rm -elementary abelian}\}$$

$$C_2:=\{ G\ |\ \exists E\ \mbox{\rm extraspecial group, } V\subset
\frac {E}{Z(E)} \  \F_p \mbox{\rm -subspace}\ |\  G\simeq  \pi^{-1}(V),$$
$${\mbox{\rm where } }\ \pi:E\to \frac {E}{Z(E)} |\
\mbox{\rm is the canonical map}\}.$$

The class $C_1$ has been described in the following result 
(cf. \cite[4.16]{Su}):

\begin{proposition} \label{suzuki}
Any $G \in C_1$ is isomorphic to one of the groups
in the following list.

a) An elementary abelian $p$-group.

b) An abelian group of type $(p,p, \dots ,p,p^2)$.

c) A central product of an extraspecial p-group $E$ and an abelian
group A. If $A$ is not elementary abelian, then $$E\cap A=Z(E)=A^p$$

\end{proposition}
 
We will need the following result from group theory. It seems to us
that it should be somewhere in the literature but we have no reference.
Although it is possible to give a proof using the classification of
extraspecial groups in section \ref{extra} and Proposition \ref{suzuki},
we give a direct proof which deals with factor systems and so is related with the algorithm 
proved for class $C_3$.  

\begin{proposition} \label{420} The 2 classes $C_i,i=1,2$ are equal.
\end{proposition}

\begin{proof}
Set $G_0=G$, $G_i=[G,G_{i-1}]$, in particular $G_1=[G,G]=G'$ and $G/N$ being 
elementary abelian implies $G_1\subset N$, so either $|G_1|=p$ or $|G_1|=1$.
As $G$ is solvable, the sequence of the $G_i$ is strictly decreasing with $G_n=\{1\}$ 
for $n>>0$. So $G_2=[G_1,G]=\{1\}$ and $G_1=G'\subset N\subset Z(G)$.
This last condition allows to define a skew-symmetric bilinear form on the $\F_p$-vector space $G/N$: if $\overline x,\overline y\in G/N$, then $[x,y]\in G'=
<\rho>\subset N$, i.e.\ $[x,y]=\rho^\epsilon$ where $\epsilon \in \Z/p\Z$. Note that $[x,y]$ is independent of the lifts of $\overline x$ and $\overline y$ to
$G$ as $N\subset Z(G)$. We define $<\overline x,\overline y>:=\epsilon $.  
Note that $x\in Z(G)$ iff $<\overline x,\overline y>=0$ for all $y\in G$.
Therefore $<.,.>$ is non degenerate iff $G'=N=Z(G)$, i.e.\ $G$ is extraspecial.

Consider the extension of groups 

$$1\to N\to G\overset{\overset{s}\curvearrowleft}{\rightarrow} V:=G/N\to 1$$

Let $s$ be a set theoretical section. To any two $v_1,v_2\in V$, we have a 
$c(v_1,v_2)\in N$ such that $s(v_1)s(v_2)=s(v_1v_2)c(v_1,v_2)$ and $c(.,.)$
is the $2$-cocycle corresponding to the equivalence class of the above 
extension in $H^2(V,N)=H^2(V,\F_p)$.
The extension is central, so $N$ has trivial action by $V$. From $c(.,.)$,
we recover $G$ in the following way: on the set $V\times \F_p$, one
defines a group structure via:
$$(v_1,\alpha)(v_2,\beta)=(v_1+v_2,\alpha+\beta +c(v_1,v_2)).$$
The form $<.,.>$ on $V$ can be expressed in terms of $c$: 
$<v_1,v_2>=[s(v_1),s(v_2)]=s(v_1)^{-1}s(v_2)^{-1}s(v_1)s(v_2)=(s(v_2)s(v_1))^{-1}s(v_1)s(v_2)=c(v_2,v_1)^{-1}s(v_2v_1)^{-1}s(v_1v_2)c(v_1,v_2)=c(v_2,v_1)^{-1}
c(v_1,v_2).$
Identifying $N$ with $\F_p$ we write $<v_1,v_2>=c(v_1,v_2)-c(v_2,v_1).$

We distinguish the cases $p>2$ and $p=2$.

Case $p>2$. 
Let $V\hookrightarrow W:=V\bigoplus V$, and $\pi:W\to V$ the projection on the first factor. We denote by $c$ a cocycle corresponding to the given group extension and 
we extend the corresponding $2$-form $<v_1,v_2>=c(v_1,v_2)-c(v_2,v_1)$ from 
$V$ to $W$ to a non degenerate skew form given by the matrix

$$
\left(\begin{array}{cccccc}
&&&-1&&\\
&A&&&\ddots&\\
&&&&&-1\\
1&&&&&\\
&\ddots&&& 0\\
&&1&&&
\end{array}\right)
$$

where $A$ is the matrix of $<.,.>$ on $V$. 

Now we obtain a $2$-cocycle $d:W\times W\to \F_p$, via 
$(w_1,w_2)\to \ <w_1,w_2>+c(\pi(w_2),\pi(w_1))$. 
We remark that $d_{|V}$ maps $(v_1,v_2)$ to $<v_1,v_2>+c(v_2,v_1)=
c(v_1,v_2)-c(v_2,v_1)+c(v_2,v_1)=c(v_1,v_2)$. So  $d_{|V}=c$ and the group
extension $E$ corresponding to $d$ therefore contains $G$ as a subgroup.

It remains to show that $|Z(E)|=p$. This amounts to the skew-form $<<.,.>>$
associated to $d$ on $W$ to be non degenerate. We compute on $W$:

$$<<w_1,w_2>>=d(w_1,w_2)-d(w_2,w_1)$$
$$=<w_1,w_2>+c(\pi(w_2),\pi(w_1))-<w_2,w_1>
-c(\pi(w_1),\pi(w_2))$$
$$=2<w_1,w_2>+<\pi(w_2),\pi(w_1)>=2<w_1,w_2>-<\pi(w_1),\pi(w_2)>.$$

Therefore $<<.,.>>$ has the matrix 

$$
\left(\begin{array}{cccccc}
&&&-2&&\\
&2A&&&\ddots&\\
&&&&&-2\\
2&&&&&\\
&\ddots&&&0\\
&&2&&&
\end{array}\right) -
\left(\begin{array}{cccccc}
&&&&&\\
&A&&&0&\\
&&&&&\\
&&&&&\\
&0&&&0&\\
&&&&&
\end{array}\right)=
\left(\begin{array}{cccccc}
&&&-2&&\\
&A&&& \ddots&\\
&&&&&-2\\
2&&&&&\\
&\ddots&&&0\\
&&2&&&
\end{array}\right)
$$

which has maximal rank as $p>2$. We conclude $E$ is extraspecial. We have obtained $G$ as subgroup of the extraspecial group $E$ and $N=Z(E)\subset G$ follows 
from the fact that in the construction above, the second factor of $V\times\F_p$
and $W\times\F_p$ correspond to $N$ and $Z(E)$ respectively.

The case $p=2$.
Using the above notation, let $n:=\dim V$ and $M_n(\F_2)$ the $\F_2$-vector 
space of $n\times n$ matrices. Any such matrix defines a bilinear form (hence
a $2$-cocycle) on $V$. Therefore we have a map of $\F_2$-vector spaces 

\begin{equation} \label{e5}
M_n(\F_2)\overset{\varphi}{\rightarrow}  H^2(V,\F_2)
\end{equation}

Moreover a matrix $A$ is in the kernel $K$ of $\varphi$ iff its associated 
$2$-cocycle $c$ defines the split extension, which is the elementary abelian 
$2$-group of rank $n+1$. 

This is equivalent to $c(v_1,v_2)=c(v_2,v_2)$ for all $v_1,v_2\in V$ and 
$c(v,v)=0$ for all $v\in V$ (here we use $p=2$).  In other words $A=A^t$
and $A$ has only zeroes on its diagonal. We conclude $\dim \ker \varphi=n-1+
n-2+...+1=\frac{n(n-1)}{2}$. It is known that $\dim
H^2(V,\F_2)=\frac{n(n-1)}{2}$ \cite{Jo} p. 169.

Therefore 
$$1\to K\to M_n(\F_2)\overset{\varphi}{\rightarrow} H^2(V,\F_2)\to 1$$ 
is exact. In particular $\varphi$ is onto, so every element of $ H^2(V,\F_2)$
can be represented by a $2$-cocycle that is a bilinear form.

Again we let the given extension correspond to the cocycle $c\in Z^2(V,\F_2)$
and by the above we may assume $c$ is bilinear corresponding to a matrix $A$.

Let $V \hookrightarrow W:=V\bigoplus V$ be the first factor and
consider $d\in Z^2(V,\F_2)$ corresponding to the matrix

$$
B=\left(\begin{array}{cccccc}
&&&1&&\\
&A&&&\ddots&\\
&&&&&1\\
&&&&&\\
&0&&& 0&\\
&&&&&
\end{array}\right)
\in M_{2n}(\F_2)
$$

and $d_{|V}=c$, so the group defined by $d$ contains $G$. 

Claim: $E$ is an extraspecial group, i.e.\ the $2$-form $<<.,.>>$ defined by 
$d$ on $W$ is non degenerate.

We calculate $<<w_1,w_2>>=d(w_1,w_2)-d(w_2,w_1)=w_1^tBw_2+w_2^tBw_1=w_1^tBw_2+
(w_2^tBw_1)^t=w_1^tBw_2+w_1^tB^tw_2=w_1^t(B+B^t)w_2$
and 

$$B+B^t=
\left(\begin{array}{cccccc}
&&1&&&\\
&A+A^t&&\ddots&\\
&&&&1&\\
1&&&&&\\
&\ddots&&&0\\
&&1&&&
\end{array}\right)
$$ 

has rank $2n$ (independently of what $A$ is).
\end{proof}

We define a third class

$$C_3:=\{ G\ |\ \exists f\in Xk[X], \ (\deg f,p)=1 \mbox{ \rm such
that } G\simeq G_{\infty,1}(f)\}$$

\begin{theorem} The 3 classes $C_i,i=1,2,3$ are equal.
\end{theorem}

\begin{proof}
It is sufficient to realize any subgroup of the extraspecial groups of type
I,II and III b. (see remark \ref{A}). So we distinguish these three
cases. The method is the following: let $E$ be an extraspecial group
and 
$G=\pi^{-1}(F\subset E/Z(E))$ a saturated subgroup. For each type of extraspecial group we consider 
a realization $C_{f_1}$ and to the sub-space $F$ we associate an additive polynomial $S(F)$ which we use in order to produce a convenient  modification of the cover   $C_{f_1}$. The key point is that our modification will not change the 
commutation rule.

Type I. So  $p>2$.  Let $E$ the extraspecial group of exponent $p$ and order $p^{2n+1}$.

Consider the realization:
$W^p-W=f_1(X):=X^{1+q}$ where $q=p^n$; then $\Ad_{f_1}(Y)=Y^{q^2}+Y$.

Note that $ E/Z(E)$ is the  automorphism group of $k[X]$ whose elements are $\sigma_y(X)=X+y$ where $y$ goes through the roots of $\Ad_{f_1}(Y)$. Then the subgroups correspond to 
those $\sigma_y$ going through $y\in F$ where $F$ is a subgroup generated by 
a subset of such roots, i.e.\ there is an additive polynomial (monic) $S_F$ which  
divides $Y^{q^2}+Y$ and $y$ goes through these roots. 
Note that necessarily $S_F$ has distinct roots so $S_F=s_0X+s_1X^p+....+
X^{p^r}$ and $s_0\neq 0$. As $\pi^{-1}(W)=E$, we can assume that $0<r<2n$.

Let us assume that $p\geq 3$ and $\el>1$ such that $(\el(\el+1),p)=1$. 

Let $f(X):=S_F(X)^{\el+1}+f_1(X)$. We remark that the conductor of 
$\red (S_F(X)^{\el+1})$  is $1+\el p^r$.  We can use the same trick as in lemma
\ref{44} for this let $j_0=1+(\el-1)p^r$. The coefficient in $\Delta(f)(X,Y)$ of 
$X^{j_0}$ is that of $\Delta(S_F^{\el+1})(X,Y)=(S_F(Y)+s_0X+...+X^{p^r})^{\el+1}-S_F(Y)^{\el+1}-S_F(X)^{\el+1}$. For this we solve the system 

$i+i_0+...+i_r=\el +1$ and $i_0+i_1p+i_2p^2+...+i_rp^r=j_0$ where $i\in\{1,2,..,\el\}$. We get $p^r-1=(i-1)p^r+i_0(p^r-1)+i_1(p^r-p)+...+i_{r-1}(p^r-p^{r-1})$. It follows that $p|i_0-1$ and $i_0\leq 1$ so $i_0=1$ and 
$i=1$, $i_1=i_2=...=i_{r-1}=0$ and so $i_r=\el-1$. Finally we have shown that 
the desired coefficient is $\frac{(\el+1)!}{1!1!(\el-1)!}S_F(Y)s_0$.

It follows that $\Ad_{f}(Y)$ divides $S_F(Y)$ which itself divides $\Ad_{f_1}(Y)=Y+Y^{q^2}$ and so  $\Ad_{f}(Y)=S_F(Y)$. 

We remark that $\Ad_{X^{\el +1}}(Y)=Y$ by lemma \ref{44}, so $\Ad_{S_F^{\el +1}(Y)}=S_F(Y)$ by \ref{49}, then by \ref{412}
for $y,z\in Z(S_F(Y))$, one has $\epsilon_f(y,z)=\epsilon_{S_F^{\el +1}}(y,z)+\epsilon_{f_1}(y,z)$. As $\epsilon_{S_F^{\el +1}}(y,z)=0$ the commutation rule is that of $E$.

The simplest choice for $p>3$ is $\el=2$ and for $p=3$,  $\el=4$.

Type II. So  $p>2$.  Let $E$ be the extraspecial group of exponent $p^2$ and order $p^{2n+1}$.
We first recall the realization of $E$ we gave in Proposition \ref{418}. Let $q=p^n$, $\theta $ a $q^2-1$-th 
root of $-1$ and  $A(X)=\theta X^q-{\theta}^qX$. Then $S(X):=A(X)+A(X)^q+...+A(X)^{q/p}$ is an additive polynomial such that $S(X)^p-S(X)=A(X)^q-A(X)=\theta^q(X^{q^2}+X)$ and  $f_0(X):=f_1(S(X))+f_2(X)$ gives a realization of $E$ where $f_1(X):=c(X^p,-X)$ and $f_2(X)=X^{1+q}$. Moreover $\Ad_{f_0}(Y)=
Y^{q^2}+Y$ so we can apply the  same strategy as for type I. 

As above we consider $\el>1$ such that $(\el(\el+1),p)=1$ and 
$f(X):=S_F(X)^{\el+1}+f_1(S(X))+f_2(X)$. If we compare to the type I case, we need to show that $\Delta (f_1\circ S)(X,Y)$ has no contribution  which cancels $\frac{(\el+1)!}{1!1!(\el-1)!}S_F(Y)s_0X^{j_0}$. 

We have ( see \ref{418}) $\Delta (f_1\circ S)(X,Y)=c(A(X)^q-A(X),A(Y)^q-A(Y))
+(F-\Id)c(S(X),S(Y))$ and $c(A(X)^q-A(X),A(Y)^q-A(Y))=c(\theta^q(X^{q^2}+X),\theta^q(Y^{q^2}+Y))=\theta^{pq}\sum_{1\leq i\leq  p-1}\frac{(-1)^{i-1}}{i}(Y^{q^2}+Y)^{p-i}(X^{q^2}+X)^i$. We remark that the equation $i_0+i_1q^2=j_0$ with 
$i_0+i_1=i$ and $1\leq i\leq p-1$ is equivalent to $i_0=1$ and $\el-1=
i_1p^{2n-r}$; so if $((\el-1)\el(\el+1),p)=1$ there is no cancellation and 
we can conclude as in the previous case. 

If $p>3$ $\el=2$ works.

If $p=3$ we need to look more carefully. In this case take $\el=4$.
Then the equation above gives $3^{2n-r}i_1=3$ which has a solution iff $r=2n-1$,
and then $i_1=1$ and $i=i_0+i_1=2$. Let us assume that $r=2n-1$. The contribution in $X^{j_0}$ is $c(Y):=2S_F(Y)s_0+ 2\theta^{3q}(Y^{q^2}+Y)$. We can write 
$Y^{q^2}+Y=(Y-\alpha)S_F(Y)$. Then $s_0\alpha =-1$ and $c(Y)=2S(Y)\theta^{3q}(Y-\alpha+\theta^{-3q}s_0)$. We remark that $\alpha^{q^2}+\alpha=0$ so 
$(-\theta^{-3q}s_0+\alpha)^{q^2}+(-\theta^{-3q}s_0+\alpha)=0$; in particular 
$c(Y)$ divides $S(Y)^2$ and we conclude as in the previous case.
 
Type III.b. So  $p=2$.  Let $E$ be the extraspecial group of type III.b. and
cardinal $p^{2n+1}$, i.e.\ it is the central product of $Q_8$ and $n-1$ copies
of $D_8$.

We have seen that the cover $W^2+W=X^{1+2^n}$ induces the extraspecial group
$E$. The corresponding additive polynomial is $Y^{q^{2}}+Y$ where $q=2^n$. 
So if $S_F(Y)$ is the additive polynomial corresponding to a saturated
subgroup we take $f(X)=S_F(X)^7+X^{1+2^n}$. We look at the contribution of 
$\Delta(S_F^{7})(X,Y)=(S_F(Y)+s_0X+...+X^{2^r})^{7}-S_F(Y)^{7}-S_F(X)^{7}$.
The contribution in $S(Y)s_0^2X^{2+2^{r+2}}$ is the only one in degree $1+2^{r+1}$
modulo multiplication by $2^\N$ and again we conclude as in the previous case.

\end{proof}

\subsection{Application}

In paragraph \ref{4.5} we haven't given a realization of extraspecial 
groups of type III.a. Such a group with cardinal
$2^{2(n-1)+1}$  is the central product $D_8*...*D_8$ ($n-1$ times). Let us
explain how we get a realization using the method above. Such a group is
a saturated subgroup of the extraspecial group
of type III.b. $Q_8*D_8*...*D_8$ of cardinal $2^{2n+1}$. The 
construction above gives the existence of $f(X)=S_F(X)^7+X^{1+2^n}$ 
where $S_F$ is an additive polynomial of degree $2^{n-1}$  such that  
the automorphism group $G_{\infty,1}(f)$
is the saturated subgroup $D_8*...*D_8$ ($n-1$ times). Note that the 
conductor is $\leq 1+6*2^{2(n-1)}$. More concretely we now give an 
explicit realization
of $D_8$ ($n=2$) with conductor 25 which is the minimal one as 
$\frac{|G_{\infty,1}(f)|}{g(f)}=\frac{2^3}{(m-1)/2}\leq \frac{2}{3}$ .

We view $D_8$ as a saturated subgroup $\pi^{-1}(F)$ of the extraspecial 
group $E:=Q_8*D_8$ where $\pi: E\to E/Z(E)=W$ for which we know 
that $f_1(X)=X^{1+2^2}$ gives a realization. The corresponding additive 
polynomial is $\Ad_{f_1}(Y)=Y^{2^4}+Y$ and $W=Z(\Ad_{f_1})\subset
\F_2^{alg}$. Then $F$ is a
subgroup of order 4 and $F=Z(S(F))$, where $S(F)$ is an additive polynomial
dividing $\Ad_{f_1}$. Therefore we can write $S(F)(Y)=Y^4+aY^2+bY$ where $b\neq 
0$.

The remainder of the division of $Y^{16}+Y$ by $Y^4+aY^2+bY$ is
$(1+b^5+ba^6)Y+(b^2a^4+ab^4+a^7)Y^2$.
Consequently we get the two equations
$1+b^5+ba^6=0 $ and
$b^2a^4+ab^4+a^7=0$.
For each couple $(a,b)$ satisfying these equations we consider $f_{a,b}:=
(X^4+aX^2+bX)^7+X^5$; then $\Ad_{f_{a,b}}(Y)=Y^4+aY^2+bY$.
Let $y\in W$ then
$P_{f_1}(X,y)+P_{f_1}(X+y,y)=y^4(y^6+y)$. We remark that $Y^6+Y$ divides 
$\Ad_{f_1}(Y)=Y^{16}+Y$ and the quotient is  $Y^{10}+Y^5+1$. This gives 
a partition
of $W$ in two sets $W_2$: The roots of $Y^6+Y$ corresponding to the $12$ 
elements of $G_{\infty,1}(f_1)$ of order $\leq 2$ and $W_4$: The roots 
of $Y^{10}+Y^5+1$ corresponding to the $20$ elements of 
$G_{\infty,1}(f_1)$ of order $4$.
Now $F:=Z(\Ad_{f_{a,b}})$ is a subgroup of $W$ and for $y\in F$ one has
$P_{f_1}(X,y)=P_{f_{a,b}}(X,y)$ and so 
$P_{f_{a,b}}(X,y)+P_{f_{a,b}}(X+y,y)=y^4(y^6+y)$. Concerning the commutation 
rule for $y,z\in F\subset W$ we have
$\epsilon_{f_{a,b}}(y,z)=\epsilon_{f_1}(y,z)=z^2y^8+zy^4+z^8y^2+z^4y=
yz(y+z)(y^2+zy+z^2)(zy^4+z^4y+1)$.

1. $a=0$ and $1+b^5=0$.
Note that in this case $\Ad_f(Y)=Y^4+bY$ and only the roots $y=0$ and 
$y=b^2$
are in $W_2$.  Moreover for $y,z\in F$ one has 
$yz(y+z)(y^2+zy+z^2)=yz(y^3+z^3)=0$; it follows that the group 
$G_{\infty,1}(f_{a,b})$ is abelian, isomorphic to
$\Z/2\Z\times \Z/4\Z$.

2. Let $A:=b^5+a^6b+1$ and $B:=b^4+a^3b^2+a^6$. Then the resultant of
$A,B$ in $b$ is $(b^5+1)(b^{10}+b^5+1)^2$. The case $b^5+1=0$ is case 1. 
above.
Now we can assume that $b^{10}+b^5+1=0$, i.e. $b$ is a primitive $15$-th 
root of 1. The equations $A=B=0$ give three sets of covers.

i) $ab=1$, i.e. $a=b^{14}$.
In this case $Y^4+aY^2+bY$ divides $Y^6+Y$, the group has exponent $2$ and
it is isomorphic to $(\Z/2\Z)^3$.

ii) $ab=b^5$ i.e. $a=b^4$.
In this case $Y^4+aY^2+bY$ has only one root ($b^7$) in common with
 $Y^{10}+Y^5+1$. It follows that the group $G_{\infty,1}(f_{a,b})$ has $2$ 
elements of order 4 so it is $D_8$. We can write $\red 
(f_{a,b})=(b^{14}+b^5)X+ 
(b+b^8)X^3+(1+b^{14}+b^{13})X^5+(b^7+1)X^7+(b^{13}+b^{10})X^9+(b^4+b+b^6)X^{11}+(b^2+1)X^{13}+b^2X^{17}+b^3X^{19}+b^9X^{21}+bX^{25} 
$, which is defined over $\F_{16}$.

iii) $ab=b^{10}$ i.e. $a=b^9$.
In this case $(Y^4+aY^2+bY)/Y$ divides $Y^{10}+Y^5+1$ and it follows that the 
group $G_{\infty,1}(f_{a,b})$ has $6$ elements of order 4, so this is $Q_8$.

\begin{remark} \rm
We could as well obtain families. For this it suffices to deal with 
$f_1$ giving a family, for example $f_1=tX^3+X^5$. The corresponding  
discussion is thought more delicate as the above.
\end{remark}

\subsection{An Algorithm} \label{algor}

Here we illustrate the algorithm which for a given $f$ gives the structure
of the group $G_{\infty,1}(f)$.
This example is a realization  $D_8$ over $\F_2$.
We have used the following Maple code:

\begin{verbatim}
 > restart;
 > f:=X^(1+2)+X^(1+2+2^2)+X^(1+2+2^4)+X^(1+2+2^5)+X^(1+2^3+2^5):
 > F:=collect(subs(X=X+Y,f),X) mod 2:
 > f1:=collect(F-subs(X=0,F)-subs(Y=0,F),[X,Y]) mod 2:
 > f2:=rem(collect(f1+f2^2-subs(X=0,f1+f2^2),[X,Y]) mod 2,X^21,X) mod 2:
\end{verbatim} 

Note that $21=40/2+1$. Here one reiterates the command until it is
stationary.

\begin{verbatim}
 > p:=f2:
 > G:=collect(f1+p^2-p,[X,Y,t]) mod 2;
G :=
(Y^24+Y^80+Y^132+Y^528+Y^192+Y^64+Y^576+Y^1280+Y^1088+Y^6+Y^3+Y^16+Y^9
+Y^272)*X^32+(Y^256+Y^128+Y^4+Y^32)*X^24+(Y^128+Y^2)*X^36+(Y^4+Y)*X^34
+(Y+Y^16)*X^40+(Y^8+Y^2)*X^33 
\end{verbatim}

Here we remark that $\Ad_f(Y)$ divides the coefficient of $X^{34}$.

\begin{verbatim}
 > G:=collect(rem(G,Y^4+Y,Y)mod 2,X);
G := 0
\end{verbatim}

Conclusion: $\Ad_f(Y)=Y^4+Y$.
 
\begin{verbatim}
> p:=collect(rem(p,Y^4+Y,Y) mod 2,X);
p :=
Y^2*X^20+Y^2*X^17+X^10*Y+Y^2*X^9+Y^3*X^8+Y^2*X^5+X^3*Y^2+Y^3*X^2+Y^2*X
 > rem(collect(p+subs(X=X+Y,p),X),Y^4+Y,Y)  mod 2;
Y^3+Y^2+Y
 > Gcd(Y^4+Y,Y^3+Y^2+Y) mod 2;
Y^3+Y^2+Y
\end{verbatim}

It follows that the $3$ roots of $Y^3+Y^2+Y$ induce $6$ order $2$ elements
and the last root $Y=1$ induces $2$ order $4$ elements.
 
\begin{verbatim}
> CY:=collect(subs(X=X+Z,p)+p,[Y]) mod 2:
 > CZ:=subs([Y=Z,Z=Y],CY):
 > C:=collect(CY-CZ,[X,Y,Z]) mod 2:
 > CC:=collect(rem(C,Y^4+Y,Y)mod 2,Z):
 > CCC:=collect(rem(CC,Z^4+Z,Z)mod 2,X);
CCC := Z^2*Y+Z*Y^2
\end{verbatim}

The group is non abelian of order $8$ with $2$ elements of order $4$; this
is $D_8$.

Note that
$f := X^3+X^7+X^{19}+X^{35}+X^{41}$ is reduced with conductor $>25$.
More generally it is a good  question to ask for realizations over $\F_2$ (i.e. $f\in \F_2[X]$)
for groups in the class $C_1$.



\begin{flushleft}
Claus LEHR\\

Michel MATIGNON \\
Laboratoire de Th\'eorie des Nombres
et d'Algorithmique Arithm\'etique,
UMR 5465 CNRS \\
Universit\'e de Bordeaux I,
351 cours de la Lib\'eration, 
33405 Talence Cedex, France \\
e-mail : {\tt matignon@math.u-bordeaux.fr}, {\tt lehr@math.u-bordeaux.fr}

\end{flushleft}


\begin{thebibliography}{Gr-Ma}
\bibitem[C-F]{C-F} R.~Carter, P.~Fong, \emph{The Sylow $2$-subgroups
of the finite classical groups}, J. of Algebra (1964), 1, 139-151.


\bibitem[De-Mu]{De-Mu} P.~Deligne, D.~Mumford, \emph{The irreducibility of the
space of curves of given genus}, Inst. Hautes Etudes Sci. Publ. Math. 
(1969), 36, 75-109.

\bibitem[De]{Deschamps} M.~Deschamps, \emph{R\'eduction
semi-stable}, Pinceaux de courbes de genre au moins deux (L.~Szpiro,
ed.), Asterisque, vol.86, (1981), pp.1-34

\bibitem[Gu]{Guralnick} R.~Guralnick, \emph{Monodromy groups of
coverings of curves}, preprint.

\bibitem[Hu]{Hu} B.~Huppert, \emph{Endliche Gruppen I}, Die Grundlehren der Mathematischen Wissenschaften, Band 134 Springer-Verlag, Berlin-New York 1967.

\bibitem[Jo]{Jo} D.L.~Johnson, \emph{Presentation of groups} Cambridge
University Press, 1976. (London mathematical society lecture note series 22). 


\bibitem[Leo]{Leo} H-W.~Leopoldt, \emph{\"Uber die Automorphismengruppe des Fermat k\"orpers}, J. Number Theory 56 (1996), no. 2, 256--282. 

\bibitem[Le-Ma1]{Le-Ma1} C.~Lehr, M.~Matignon, \emph{Wild monodromy
and automorphisms of curves}, conference proceedings, Tokyo, (2002)
(T.~Sekiguchi, N.~Suwa, editors), to appear.

\bibitem[Le-Ma2]{Le-Ma2} C.~Lehr, M.~Matignon, \emph{Wild monodromy and automorphisms of curves}, in preparation.

\bibitem[Le1]{Le1} C.~Lehr, \emph{ Reduction of $p$-cyclic Covers of the Projective Line},  Manuscripta Math. 106 (2001) 2, 151-175.

 
\bibitem[Liu]{Liu} Q.~Liu, \emph{ Algebraic Geometry and Arithmetic Curves}, Oxford Graduate Texts in 
Mathematics, 6 (2002), Oxford University Press.
  
\bibitem[Ma]{Ma} M.~Matignon, \emph{ Vers un algorithme pour la
r\'eduction stable des rev\^etements $p$-cycliques de la droite
projective sur un corps $p$-adique}, 
Mathematische Annalen 325, 323-354 (2003).

\bibitem[Po1]{Poonen1} B.~Poonen, \emph{Varieties without extra
automorphisms I: Curves}, Math. Res. Lett. 7 (2000), no. 1, 67--76.

\bibitem[Po2]{Poonen2} B.~Poonen, \emph{Varieties without extra
automorphisms II: Hyperelliptic curves}, Math. Res. Lett. 7 (2000), 
no. 1, 77--82. 

\bibitem[St1]{St1} H.~Stichtenoth, \emph{\"Uber die Automorphismengruppe eines algebraischen Funktionenk\"orpers von Primzahlcharakteristik. I. Eine Absch\"atzung der Ordnung der Automorphismengruppe.} Arch. Math. 24 (1973) 527--544.

\bibitem[St2]{St2} H.~Stichtenoth, \emph{\"Uber die Automorphismengruppe eines algebraischen Funktionenk\"orpers von Primzahlcharakteristik. II. Ein spezieller Typ von Funktionenk\"orpern.} Arch. Math. 24 (1973), 615--631.

\bibitem[Su]{Su} M.~Suzuki, \emph{Group theory II}, Grundlehren der Mathematischen Wissenschaften 248. Springer-Verlag, New York, 1986.

\bibitem[vdG-vdV]{vdG-vdV}
G.~van der Geer, M.~van der Vlugt, \emph{Reed-Muller 
codes and supersingular curves I.} Compositio Math. (1992), 84,
no. 3, 333--367.                                              


\end{thebibliography}
\end{document}